\documentclass{amsart}
\usepackage{graphicx,color}
\newtheorem{thm}{Theorem}[section]
\newtheorem{cor}[thm]{Corollary}
\newtheorem{lem}[thm]{Lemma}
\newtheorem{prop}[thm]{Proposition}
\theoremstyle{definition}
\newtheorem{defn}[thm]{Definition}
\theoremstyle{remark}
\newtheorem{rem}[thm]{Remark}
\numberwithin{equation}{section}

\newcommand{\abs}[1]{\left\vert#1\right\vert}

\newcommand{\eps}{\varepsilon}

\newcommand{\A}{\mathcal{A}}
\newcommand{\tu}{\tilde{u}}

\newcommand{\wind}{{\rm wind}}
\newcommand{\ind}{{\rm ind}}
\newcommand{\C}{\mathbb{C}}
\newcommand{\R}{\mathbb{R}}
\newcommand{\N}{\mathbb{N}}
\newcommand{\Z}{\mathbb{Z}}
\newcommand{\pr}{{\rm pr}}
\newcommand{\w}{{\rm w}}
\newcommand{\gr}{{\rm gr}}
\newcommand{\Hom}{{\rm Hom}}
\newcommand{\Giroux}{{\rm Giroux}}
\newcommand{\FF}{\mathcal{F}}
\newcommand{\MM}{\mathcal{M}}

\newcommand{\p}{\partial}
\newcommand{\om}{\omega}
\newcommand{\Om}{\Omega}

\begin{document}
\parindent=0pt
\parskip=2pt

\title[Planar Weinstein Conjecture]{The Weinstein Conjecture for Planar
Contact Structures in Dimension Three}%
\author{Casim Abbas$^1$}%
\address{}%
\email{}%
\author{Kai Cieliebak}%
\address{}%
\email{}%
\author{Helmut Hofer$^2$}%
\address{}%
\email{}%
\thanks{$^1$ Research partially supported by a Michigan State University IRGP grant}
\thanks{$^2$ Research partially supported by NSF Grant DMS-0102298}


\date{March 4, 2005}
\begin{abstract}
\parindent=0pt
\parskip=2pt
In this paper we describe a general strategy
for approaching
the Weinstein conjecture in dimension three.
We apply this
approach to prove the Weinstein conjecture for a new class of
contact manifolds (planar contact manifolds). 
We also discuss how the present approach reduces the
general Weinstein conjecture in dimension three to a compactness
problem for the solution set of a first order elliptic PDE.


\end{abstract}
\maketitle
\tableofcontents

\section{Introduction}\label{sec:intro}
The following considerations are part of the program initiated in
\cite{Hofer-Weinstein-conj} and extended in \cite{Hofer2000} of proving
the general Weinstein conjecture in dimension three. The key
observation in \cite{Hofer-Weinstein-conj} was the equivalence
between the assertion of the Weinstein conjecture and the
existence of a non-constant holomorphic curve for a suitable
nonlinear Cauchy-Riemann type equation. As discussed in
\cite{Hofer2000}, this equivalence has its limitations. However, it was
suggested that a suitable modification of the
holomorphic curve equation should be the key to a
proof of the general Weinstein conjecture in dimension $3$. In the
current joint work of the authors, the proof of the general
Weinstein conjecture in dimension three has been reduced to a
compactness question of certain moduli spaces for the generalized
holomorphic curve equation. As we know from Giroux's work,
any (co-oriented) contact structure is supported by an open book
decomposition. In our approach, the compactness problems only arise
if the pages of the open book decomposition are non-planar (i.e., of
positive genus). If the pages are planar these difficulties do not
arise. In this paper we describe our approach for this particular
case.

\subsection{Versions of the Weinstein conjecture}
Before we give more details we start by providing the necessary
background. Consider a closed three-manifold $M$ equipped
with a contact structure $\xi$. In this paper we assume all contact
structures to be cooriented, i.e., $\xi=\ker(\lambda)$ is defined by a
contact 1-form $\lambda$. We denote the associated Reeb vector field
by $X_{\lambda}$. Recall that the (generalized) three-dimensional
Weinstein conjecture states the following, see \cite{Weinstein}:

\noindent {\bf Conjecture (A.~Weinstein, 1978). }
Every Reeb vector field $X$ on a
closed three-dimensional manifold $M$ admits a periodic orbit.

In fact, Weinstein added the hypothesis that the first
cohomology group $H^1(M;{\mathbb R})$
vanishes, but there is no indication that this additional
hypothesis is needed. Moreover, Weinstein made his conjecture for
Reeb vector fields on odd-dimensional manifolds of arbitrary
dimensions. We point out at there are strong indications
that in fact a stronger form of the Weinstein conjecture is true,
which we again formulate in the three-dimensional case:

\noindent {\bf Strong version of the Weinstein conjecture. }
For every Reeb vector field $X$ on a closed three-dimensional
manifold $M$ there exist finitely many periodic orbits
$(x_i,T_i)$, $i=1,..,n$, so that the first homology classes
$[x_1],...,[x_n]$ induced by the loops $x_i:{\mathbb
R}/(T_i{\mathbb Z})\rightarrow M$ sum up to $0$:
$$
\sum_{i=1}^{n} [x_i] =0.
$$
Here the periods $T_i>0$ need not to be the minimal periods.

We will say that the {\em (strong) Weinstein conjecture holds for a
contact form $\lambda$} if the associated Reeb vector field
satisfies the conclusion of the (strong) Weinstein conjecture.

\subsection{Generalized holomorphic curve equations}
We write $\pi:TM\rightarrow \xi$ for the projection along
$X_{\lambda}$.
Fix a complex structure $J$ on $\xi$ such that
$d\lambda(\cdot,J\cdot)$
defines a positive definite metric on $\xi$. We will call
such complex structures {\em compatible (with $d\lambda$)}. Let
us begin with an assertion  reducing the Weinstein conjecture to the
study of the following nonlinear first order elliptic system. The
solutions of interest are 5-tuples
$(S,j,\Gamma,\tilde{u},\gamma)$ consisting of a closed Riemann
surface $(S,j)$, a finite subset $\Gamma\subset S$, a proper map
$\tilde{u}=(a,u):\dot{S}\rightarrow {\mathbb R}\times M$, where
$\dot{S}=S\setminus\Gamma$, and a one-form $\gamma $ on $S$ so
that
\begin{equation}\label{eq1}
\left\{
\begin{array}{ccc}
&\pi\circ Tu\circ j = J\circ \pi\circ Tu\ \hbox{on}\ \dot{S},&\\
&(u^{\ast}\lambda)\circ j = da +\gamma\ \hbox{on}\ \dot{S},&\\
& d\gamma= d(\gamma\circ j)=0\ \hbox{on}\ S,&\\
&E(\tilde{u})<\infty.&
\end{array}
\right.
\end{equation}
Here the {\em energy} $E(\tilde{u})$ is defined by
\begin{eqnarray*}
E(\tilde{u})=\sup_{\varphi\in\Sigma}\ \int_{\dot{S}}\
\tilde{u}^{\ast}d(\varphi\lambda),
\end{eqnarray*}
where $\Sigma$ consists of all smooth maps $\varphi:{\mathbb
R}\rightarrow [0,1]$ with $\varphi'(s)\geq 0$ for all $s\in
{\mathbb R}$.

The following theorem, which is an easy modification of a result
by Hofer~\cite{Hofer-Weinstein-conj,Hofer2000}, shows that the
Weinstein conjecture is equivalent to an existence result for a
generalized holomorphic curve (we restrict ourselves to the case of
three dimensions in the following discussion):
\begin{thm}
Let $(M,\lambda)$ be a closed three-dimensional manifold equipped
with a contact form $\lambda$. Then the associated Reeb vector
field has periodic orbits if and only if the associated PDE-problem
(\ref{eq1}) has a non-constant solution.
\end{thm}

Note, however, that a nontrivial solution need not to have any
puncture due to the harmonic perturbation $\gamma$. If the Riemann
surface $(S,j)$ occurring in (\ref{eq1}) is a sphere
it follows immediately that $\gamma=0$, and the generalized equation
reduces to the usual equation for punctured holomorphic curve
in symplectizations.

\subsection{Open book decompositions and the main result}
An {\em open book decomposition} of a closed 3-manifold $M$ is a pair
$(L,\pr)$ consisting of a fibered link $L\subset M$ (the {\em
  binding}) and a fibration
$\pr:M\setminus L\rightarrow S^1$ whose fibers $\pr^{-1}(t)$ (the {\em
  pages}) are the interiors of smooth compact embedded surfaces
in $M$ bounded by $L$.

\begin{defn}\label{qqq1}
Following~\cite{Giroux-lecture}, we say that a contact structure $\xi$
on a  closed 3-manifold $M$ is {\em supported by an open book
decomposition}  $(L,\pr)$ if there exists a contact form $\lambda$
defining $\xi$ so that:
\begin{itemize}
\item The form $d\lambda$ induces an area form on each leaf $F$ of $\pr$.
\item The form $\lambda$ defines a volume form on $L$ inducing
 the orientation as boundary of $(F,d\lambda)$.
\end{itemize}
We will call $\lambda$ a {\em Giroux form} associated to $(L,\pr)$ and
denote such forms by $\lambda_{Giroux}$.
\end{defn}

\begin{rem}
(1) The definition implies that each component of $L$ is a periodic orbit
of the Reeb vector field associated to the Giroux form. We call the
components of $L$ the {\em binding orbits}. Note that any Reeb orbit
which is not a binding orbit hits any page in forward and backward
time. 

(2) A given Giroux form can be modified
near the binding $L$ to have additional properties.
For example, we can arrange that each binding orbit has a
neighborhood isomorphic to that of a periodic orbit in the round
sphere $S^3$. Alternatively, we can arrange that the binding
orbits are nondegenerate elliptic periodic orbits
(see~\cite{Abbas}). 

(3) Multiplying a Giroux form by some positive number we obtain
another Giroux form.
\end{rem}

Giroux's fundamental result is the following,
see~\cite{Giroux-lecture,Giroux-ICM}.

\begin{thm}
Any (co-orientable) contact structure on a closed 3-manifold $M$ is
supported by an open book decomposition.
\end{thm}

\begin{defn}
Let us call an open book decomposition {\em planar} if its pages have
genus zero. Call a contact structure {\em planar} if it is supported
by a planar open book decomposition.
\end{defn}

The main result in this paper is the following:

\begin{thm}[Strong Weinstein Conjecture for Planar Contact
  Structures]\label{thm:main}
Let $\xi$ be a planar contact structure on an oriented closed
three-manifold $M$. Then the strong version of the Weinstein
conjecture holds for any contact form defining $\xi$.
\end{thm}

In view of this theorem and the results
of~\cite{Hofer-Weinstein-conj}, the Weinstein conjecture is now
established for every contact form defining a contact
structure $\xi$ on a closed oriented 3-manifold $M$ if at least one of
the following conditions is met:
\begin{itemize}
\item[(1)] The contact structure $\xi$ is overtwisted
  (\cite{Hofer-Weinstein-conj}).
\item[(2)] The second homotopy group of $M$ is nontrivial
  (\cite{Hofer-Weinstein-conj}).
\item[(3)] The contact structure $\xi$ is planar (present paper).
\end{itemize}

\begin{rem}
(1) Recent progress in the understanding of contact three-manifolds,
most notably an important result by Eliashberg~\cite{Yasha}, has led
to serious advances in the study of the "planarity question". Indeed,
in a recent paper ~\cite{Etnyre} J.~Etnyre shows that not all contact
structures are planar. He also shows that every overtwisted contact
structure is planar, so case (1) above is a consequence of case
(3). The Weinstein conjecture remains open for tight contact forms on
closed 3-manifolds with vanishing second homotopy group for which the
underlying contact structure is not planar.

(2) It was pointed out by J.~Etnyre that one can modify our proof in
the planar case by putting on top of our construction Eliasherg's
symplectic cobordism~\cite{Yasha}. Then one can work
with honest spheres rather than punctured spheres. The proof
then has to make use of positivity of intersections, adjunction
formula, self-intersection index, automatic transversality, and
the compactness results for punctured holomorphic curves
in~\cite{BEHWZ}. Our
arguments may be viewed as relative versions of these concepts. 

The use of Eliashberg's cobordism would somewhat simplify the
arguments in the planar case if one does not like to work with
non-compact curves. However, an
index calculation shows that Eliashberg's construction does not
help to prove the Weinstein conjecture in non-planar cases,
whereas our constructions are designed precisely for this case.
The only problem at the moment is the lack of a compactness proof for
the generalized equation~(\ref{eq1}), although
we are making progress on this question.
%
\end{rem}

{\bf Acknowledgement:} The third author would like
to thank Richard Siefring for helpful discussions on intersection
questions which simplified some of our arguments.

\section{Recollections on finite energy
  spheres}\label{sec:finite-energy-spheres}

In this section we collect some facts about solutions of the
PDE~\eqref{eq1}. Most of the results needed are scattered in the
literature. Some of them need additional explanations and are further
discussed.
Throughout, $M$ is a closed oriented 3-manifold, $\xi$ is a contact
structure defined by a contact form $\lambda$, $J$ is a compatible
complex multiplication on $\xi$, and $\pi:TM\to\xi$ is the projection
along the Reeb vector field $X_\lambda$. {\em From now on, we will
restrict ourselves to planar curves}, i.e., the surface $S$
in~(\ref{eq1}) is diffeomorphic to the sphere. Then the harmonic form
$\gamma$ vanishes and the PDE~(\ref{eq1}) reduces to
\begin{equation}\label{eq2}
\left\{
\begin{array}{ccc}
&\pi\circ Tu\circ j = J\circ \pi\circ Tu,&\\
&(u^{\ast}\lambda)\circ j = da,&\\
&E(\tilde{u})<\infty.&
\end{array}
\right.
\end{equation}
A solution $(S,j,\Gamma,\tilde u)$ of equation~(\ref{eq2}) is called a
{\em (special) finite energy sphere}. Equation~\eqref{eq2} can be
written in a more concise form as follows. Associate to $J$ the
almost complex structure $\tilde J$ on $\R\times M$ defined by
$$
   \tilde J|_\xi := J:\xi\to\xi,\qquad \tilde J\frac{\p}{\p
   r} := X_\lambda,\qquad \tilde J X_\lambda := -\frac{\p}{\p r},
$$
where $r$ denotes the coordinate on $\R$. Note that $\tilde J$ is
$\R$-invariant and compatible with the symplectic form $d(e^r\lambda)$
in the sense that $d(r^r\lambda)(\cdot,\tilde J\cdot)$ defines a
Riemannian metric. Then equation~\eqref{eq2} is equivalent to
$$
   T\tilde u\circ j=\tilde J\circ T\tilde u,\qquad E(\tilde
   u)<\infty.
$$
We also need to consider a generalization of
equation~\eqref{eq2}. Let $\lambda^+,\lambda^-$ be two contact forms
defining the same contact structure $\xi$ such that
$$
   \lambda^+ = f^+\cdot\lambda^-
$$
for a function $f^+>1$ on $M$. Pick a positive function $f$ on
$\R\times M$ and a constant $R>0$ such that $\frac{\p f}{\p r}\geq 0$
and
$$
   f(r,x) = \begin{cases}
      f^+(x) & \text{ for }r\geq R, \cr
      1 & \text{ for }r\leq -R.
   \end{cases}
$$
Note that $\om_f:=d(e^rf\lambda^-)$ is a symplectic form on $\R\times
M$. Let $\tilde J$ be a {\em compatible
almost complex structure} on $\R\times M$. This means that
$\omega_f(\cdot,\tilde J\cdot)$ defines a Riemannian metric. Moreover,
we assume that
$$
   \tilde J = \begin{cases}
      \tilde J^+ & \text{ on }[R,\infty)\times M, \cr
      \tilde J^- & \text{ on }(-\infty,-R]\times M,
   \end{cases}
$$
where $\tilde J^\pm$ are the $\R$-invariant almost complex structures
associated to complex multiplications $J^\pm:\xi\rightarrow\xi$
compatible with $\lambda^\pm$. We now study smooth maps
$\tilde{u}:\dot{S}\rightarrow {\mathbb R}\times M$
satisfying
\begin{equation}\label{eq3}
   T\tilde{u}\circ j = \tilde{J}\circ T\tilde{u}, \qquad E(\tilde
   u)<\infty.
\end{equation}
Here the energy $E(\tilde{u})$ is defined by
\begin{eqnarray*}
E(\tilde{u})=\sup_{\varphi\in\Sigma}\ \int_{\dot{S}}\
\tilde{u}^{\ast}d(\varphi f\lambda^-),
\end{eqnarray*}
with $\Sigma$ as before. We call solutions of equation~(\ref{eq3}) {\em
generalized finite energy spheres}.

\subsection{The Reeb flow near a periodic orbit}\label{ss:Reeb-flow}
Let $x$ be a periodic Reeb orbit on $(M,\lambda)$ of period $T$.
Denote by $\phi_t:M\to M$ the Reeb flow, thus
$\phi_T\bigl(x(0)\bigr)=x(0)$. The linearized Reeb flow along $x$
gives rise to a family of linear maps $\Phi_t:\xi_{x(0)}\to\xi_{x(t)}$
which preserve the symplectic form $d\lambda|_\xi$. We call $x$ {\em
  nondegenerate} if $\Phi_T:\xi_{x(0)}\to\xi_{x(0)}$ does not have $1$
in its spectrum. Then two cases can occur. Either both eigenvalues are
real, then we call $x$ {\em hyperbolic}, or both are non-real, then
$x$ is called {\em elliptic}.

Closely related to the linearized Reeb flow is the {\em asymptotic
operator along $x$}
\begin{equation}\label{asy-operator}
   A\eta := -J(x)(\nabla_t\eta-\nabla_{\eta}X_\lambda)
\end{equation}
acting on sections $\eta(t)=\eta(t+T)$ of the bundle $x^*\xi$. Its
kernel corresponds to eigenvectors of $\Phi_T$ with eigenvalue $1$,
so in the nondegenerate case the kernel is trivial. Moreover,
eigenfields of $A$ have no zeroes. Fix a
trivialization of the bundle $x^*\xi$. In this trivialization
each eigenfield of the self-adjoint operator $A$ has a winding
number which depends only on the eigenvalue, see~\cite{Hofer-Kriener}
for details. The winding number
increases with the eigenvalue, and each winding number occurs for
precisely two eigenvalues (counted with multiplicities). Denote by
$\alpha(x)$ the winding number corresponding to the largest
negative eigenvalue. If $x$ is nondegenerate and elliptic both
eigenvalues with winding number $\alpha(x)$ are negative and the
{\em Conley-Zehnder index} of $x$ (in the given trivialization,
see~\cite{Hofer-Kriener}) is given by
$$
    \mu(x) = 2\alpha(x) + 1.
$$
We will also need a weighted version of this relation. For a {\em
weight} $\delta<0$ which is not an eigenvalue define the {\em weighted
Conley-Zehnder index} by
$$
    \mu_\w(x) := 2\alpha_\w(x) + 1,
$$
where $\alpha_\w(x)$ is the winding number corresponding to the
largest eigenvalue $<\delta$. Note that $\alpha_w(x)$ and $\mu_\w(x)$
are the winding number corresponding to the largest negative
eigenvalue, respectively Conley-Zehnder index, of the {\em weighted
asymptotic operator} $A_\w:=A-\delta$.

\subsection{Asymptotics near a puncture}\label{ss:asymptotics}
Next we describe the behavior of solutions of~(\ref{eq2}) near a
positive puncture. The reference for this section
is~\cite{HWZ-asymptotics}. Let $x$ be a periodic Reeb orbit of period
$T>0$. Denote by $\tau>0$ the
minimal period and by $k\in\N$ its covering number, so that
$T=k\tau$. In suitable local coordinates in a
tubular neighborhood $U$ of $x$ the contact form is given by
$$
\lambda= f(d\vartheta +x\,dy)
$$
where $(\vartheta,x,y)\in S^1\times\R^2$ with $S^1=\R/\Z$. Here the
periodic orbit $x$ corresponds to $t\mapsto
(kt,0,0)\in S^1\times\R^2$, and $f>0$ is a function satisfying
$f(\vartheta,0,0)=\tau$ and $df(\vartheta,0,0)=0$.

Let $\tilde{u}=(a,u):[0,\infty)\times
S^1\rightarrow {\mathbb R}\times M$ be a solution of~(\ref{eq2}) such
that $u(s,\cdot)\rightarrow x$ and $a(s,\cdot)\to+\infty$ as
$s\to\infty$. After replacing 
$[0,\infty)\times S^1$ by $[R,\infty)\times S^1$ for a sufficiently
large $R$, we may assume that the image of $u$ is contained in a
neighbourhood $U\subset S^1\times {\mathbb R}^2$ above. Hence we can
write
$$
\tilde{u}(s,t)=(a(s,t),\vartheta(s,t),z(s,t))
$$
in the coordinates above, with $z=(x,y)\in\R^2$. The following
asymptotic behaviour was established in~\cite{HWZ-asymptotics}.

\begin{thm}\label{thm:asymptotics}
Suppose that $x$ is nondegenerate of period $T$ and covering number
$k$. Then there exist constants
$a_0, \vartheta_0\in\R$ and $d > 0$ such that
\begin{align*}
  \abs{\partial^{\beta} [a(s,t) - Ts - a_0 ]} & \leq C_\beta e^{-ds}\\
  \abs{\partial^{\beta} [\vartheta
   (s,t)-kt -\vartheta_0 ] }   & \leq C_\beta e^{-ds}
\end{align*}
for all multi-indices $\beta$, with constants  $C_\beta$ depending on
$\beta$.
Moreover, if the $z$-part does not vanish identically we have the
asymptotic formula for the transversal approach to $x(t)$:
$$
  z(s,t) = e^{\int^s_{s_0} \lambda(\sigma) d \sigma} [e(t) + r (s,t)] \in
  \R^2,
$$
where $\partial^{\beta} r(s,t) \rightarrow 0$ as $s \rightarrow
\infty$, uniformly in $t$ for all derivatives. Here $\lambda:
[s_0, \infty) \rightarrow \R$ is a smooth function
satisfying
$$
  \lambda(s) \rightarrow \lambda < 0 \qquad \mbox{as $s \rightarrow +
  \infty$},
$$
where $\lambda<0$ is an eigenvalue of the asymptotic operator $A$
along $x$ defined in the previous section and $e(t) = e(t+1) \neq 0$
is an eigenfield to $\lambda$.
\end{thm}

In particular, this implies that $u:\dot S\to M$ admits a continuous
extension $\bar{u}$ to the circle compactification $\bar{S}$  of
its domain $\dot{S}$. The behaviour near a negative puncture (at which
$a\to-\infty$) is similar. In the following we will only need positive
punctures. Note that equation~\eqref{eq3} agrees with
equation~\eqref{eq2} for $|r|\geq R$, so Theorem~\ref{thm:asymptotics}
also applies to generalized finite energy spheres.

\subsection{Linear Fredholm theory}\label{ss:lin-Fredholm}
Following~\cite{HWZ-Fredholm}, we introduce a special class of linear
Fredholm operators over a punctured Riemann sphere $(\dot S,j)$. They
act on sections of a trivial complex line bundle $V=\dot S\times\R^2$
with fibrewise complex structure $i=i(z)$, $z\in\dot S$. Denote by
$\A_0\to\dot S$ the bundle of complex antilinear bundle homomorphisms
$T\dot S\to V$. Let $C$ be a smooth section in the
bundle $\Hom_\R(V,\A_0)\to\dot S$ of real bundle homomorphisms
$V\to\A_0$. We call $C$ {\em admissible} if at every
puncture the following holds. (We assume all punctures to be positive,
although for this subsection this makes no difference). Let $(s,t)$ be
polar coordinates such that $s\to\infty$ at the puncture. Then there
exist smooth loops of complex structures $i^+(t)$ on $\R^2$ and
$2\times 2$-matrices $C^+(t)$ such that
$$
   i(s,t)\to i^+(t),\qquad C(s,t)\cdot\frac{\p}{\p s}\to C^+(t)
$$
in $C^\infty$ as $s\to\infty$. The matrices $C^+(t)$ are symmetric
with respect to the metrics $\om\bigl(\cdot,i^+(t)\cdot\bigr)$, where
$\om$ is the standard symplectic form on $\R^2$. Moreover, we require
that the asymptotic operator
$$
   A^+\eta := -i^+(t)\frac{\p\eta}{\p t}-C^+(t)\eta
$$
acting on smooth functions $\eta:S^1\to\R^2$ has trivial kernel. Thus
the equation $A^+\eta=0$ defines a path of symplectic $2\times
2$-matrices $\Phi_t$ such that $\Phi_0=\mbox{Id}$ and $\Phi_1$ does not
have $1$ in its spectrum. Denote by $\mu^+$ the Conley-Zehnder index
of this path.

We associate to an admissible $C$ the operator
$L_C:\Om^0(V)\to\Om^0(\A_0)$
$$
   L_Cv := Tv + i\circ Tv\circ j + Cv
$$
acting on sections of the bundle $V\to\dot S$. Let $\#\Gamma$ be the
number of (positive) punctures of $\dot S$ and
$$
   \mu(L_C) := \sum_j\mu_j^+
$$
the sum of the Conley-Zehnder indices at the punctures. The following
result was proved in~\cite{HWZ-Fredholm}.

\begin{prop}\label{prop:lin-Fredholm}
The operator $L_C$ associated to an admissible $C$ defines a Fredholm
operator $L_C:E\to F$ between suitable Sobolev (or H\"older)
completions of $\Omega^0(V)$ and $\Omega^0({\mathcal A}_0)$ of index
$$
   \ind(L_C) = \mu(L_C) + 2 - \#\Gamma.
$$
\end{prop}

The arguments in~\cite{HWZ-asymptotics} show that elements in the
kernel of $L_C$ have asymptotics at a puncture analogous to the
$\xi$-component $z$ in Theorem~\ref{thm:asymptotics}.

\begin{cor}\label{cor:asymptotics}
A nontrivial element $v$ in the kernel of $L_C$ has the following
asymptotic behaviour in polar coordinates near a puncture:
$$
  v(s,t) = e^{\int^s_{s_0} \lambda(\sigma) d\sigma} [e(t) + r (s,t)]
  \in \R^2,
$$
where $\partial^{\beta} r(s,t) \rightarrow 0$ as $s \rightarrow
\infty$, uniformly in $t$ for all derivatives. Here $\lambda:
[s_0, \infty) \rightarrow \R$ is a smooth function satisfying
$$
  \lambda (s) \rightarrow \lambda < 0 \qquad \mbox{as $s \rightarrow +
  \infty$},
$$
where $\lambda<0$ is an eigenvalue of the asymptotic operator $A^+$ at
the puncture and $e(t) = e(t+1) \neq 0$ is an eigenfield to $\lambda$.
\end{cor}

Let us discuss the effect of exponential weights. By
Corollary~\ref{cor:asymptotics}, a nontrivial element in the kernel of
$L_C$ approaches zero at the $j$-th puncture
with an exponential rate given by an eigenvalue $\lambda_j<0$ of the
asymptotic operator $A_j^+$. For weights $\lambda_j<\delta_j<0$ that
are not eigenvalues, denote by $E_\w,F_\w$ the weighted Sobolev spaces
of sections converging to zero at the punctures with exponential rates
$\delta_j$ or faster. Thus
$E_\w$ is the space of sections $\eta$ in $V\to\dot S$
such that $\eta_\w\in E$, where $\eta_\w$ is defined by
multiplying $\eta$ by a positive smooth
function which agrees with $e^{-\delta_js}$ near the $j$-th puncture,
and $F_\w$ similarly. Define the {\em weighted Fredholm index}
$\ind_\w(L_C)$ as the index of the linear Fredholm operator
$L_C:E_\w\to F_\w$. Note that $\eta\to\eta_\w$ defines an isomorphism
$E_\w\to E$ (and similarly for $F$) which conjugates the operator
$L_C:E_\w\to F_\w$ to the operator $L_{C_\w}:E\to F$ associated to an
admissible $C_\w$. A
simple computation shows that the asymptotic operators of $C_\w$ are
precisely the weighted asymptotic operators
(cf.~Section~\ref{ss:Reeb-flow}) at the punctures. Hence by
Proposition~\ref{prop:lin-Fredholm},
\begin{equation}\label{eq:weighted-ind}
   \ind_\w(L_C) = \ind(L_{C_\w}) = \mu_\w(L_C) + 2 - \#\Gamma,
\end{equation}
where $\mu_\w(L_C)$ is the sum of the weighted Conley-Zehnder indices
at the punctures.

Next consider a nontrivial element $v$ in the kernel of $L_C$.
By Corollary~\ref{cor:asymptotics}, it converges to zero at the $j$-th
puncture from the direction of some eigenvector $e_j$ of the
asymptotic operator $A_j^+$. Denote the
winding number of $e_j$ by $w_j^+$ and define the {\em
winding number} of $v$ by
$$
    \wind(v) := \sum w_j^+.
$$
The asymptotics of $v$ and the similarity principle imply
(cf.~\cite{HWZ-Fredholm})

\begin{lem}\label{lem:lin-winding}
Let $v$ be a nontrivial element in the kernel of $L_C$. Then $v$ has
only finitely many zeroes, each zero has positive multiplicity, and
their algebraic sum equals $\wind(v)$.
\end{lem}


\subsection{Nonlinear Fredholm theory}\label{ss:nonlin-Fredholm}
Next we recall the Fredholm theory for equation~\eqref{eq3}. The basic
references are~\cite{HWZ-Fredholm} for the embedded case (which is all
we need), and~\cite{Dragnev} for the general case. Let $\tilde
u=(a,u):\dot S\to R\times M$ be a generalized finite energy sphere with
asymptotic orbits $x_j$. {\em From now on we assume that all the
punctures are positive (i.e., $a\to+\infty$) and all the asymptotic
orbits are distinct, simple and nondegenerate elliptic.} Denote by
$\mu_j$ and $\alpha_j$ 
their Conley-Zehnder indices, respectively winding numbers of the
largest negative eigenvalue, with respect to trivializations induced
by a trivialization of $u^*\xi$ over $\dot S$. So we have
$$
    \mu_j = 2\alpha_j + 1.
$$
Denote by $\mu(\tilde u):=\sum\mu_j$
the {\em Conley-Zehnder index of $\tilde u$} and by $\#\Gamma$ the
number of (positive) punctures of $\dot S$.

Denote by $\MM$ the moduli space of solutions of equation~\eqref{eq3}
with $\#\Gamma$ positive punctures and asymptotic orbits $x_j$.
The space
$\MM$ can be described as the zero set of the nonlinear Cauchy-Riemann
operator defined by~\eqref{eq3} on a suitable Banach manifold of maps
$\dot S\to\R\times M$ times the moduli space $\MM_{0,\#\Gamma}$ of
$\#\Gamma$ points on the sphere. Its linearization at $\tilde u$ is a linear
Fredholm operator $D:E\times T\MM_{0,\#\Gamma}\to F$ between
Banach spaces. Here $E$ and $F$ are suitable
Sobolev completions of the space of sections, respectively
$(0,1)$-forms, in the pullback bundle $\tilde u^*T(\R\times
M)$. According to~\cite{HWZ-Fredholm} for embeddings $\tilde u$,
and~\cite{Dragnev} in general, the Fredholm index of $D$ is
given by
\begin{equation}\label{eq:full-ind}
    \ind(\tilde u) = \mu(\tilde u) - 2 + \#\Gamma.
\end{equation}

If $\tilde u$ is an embedding there is an alternative
description developed in~\cite{HWZ-Fredholm}. Write nearby curves as
graphs of sections in the complex normal bundle $N\to\dot S$ to
$C=\tilde u(\dot S)$ in ${\mathbb R}\times M$. Equation~\eqref{eq3}
translates into a Monge-Ampere type equation for sections of $N$ whose
linearization $D^N$ at the zero-section is the projection
of $D$ onto $N$.
Set $\mu^N(\tilde u):=\sum\mu_j^N$, where the {\em normal
  Conley-Zehnder indices} 
$$
   \mu_j^N=2\alpha_j^N+1
$$ 
are computed with
respect to trivializations induced by a trivialization of the normal
bundle $N\to\dot S$. Note that in view of the asymptotics
(Theorem~\ref{thm:asymptotics}) the bundles $u^*\xi$ and $N$ agree
near the punctures. Comparing these bundles over $\dot S$ yields the
following relation between the Conley-Zehnder index and the normal
Conley-Zehnder index, see~\cite{HWZ-asymptotics}:
\begin{equation}\label{eq:ind-normal-ind}
   \mu(\tilde{u})=\mu^N(\tilde{u}) + 4 -2\sharp\Gamma.
\end{equation}
It allows us to express the index of $\tilde u$ in terms of the normal
Conley-Zehnder index:
\begin{equation}\label{eq:normal-ind}
   \ind(\tilde u) = \mu^N(\tilde u) + 2 - \#\Gamma.
\end{equation}
On the other hand, the operator $D^N$ is an admissible operator of the
form considered in Section~\ref{ss:lin-Fredholm} (see the proof of
Lemma~\ref{lem:nearby-sol} below), and by
Proposition~\ref{prop:lin-Fredholm} its index is given by the
right-hand side of equation~\eqref{eq:normal-ind}.
This must of course be the case because $D$ and $D^N$ both
describe nearby solutions of the same equation~\eqref{eq3}.
The following lemma is implicit in~\cite{HWZ-Fredholm}.

\begin{lem}\label{lem:nearby-sol}
Let $\tilde u$ be an embedded solution of~\eqref{eq3} with only
positive punctures asymptotic to elliptic Reeb orbits $x_j$.
Let $v$ be a section in the normal bundle $N\to\dot S$ whose graph
describes a solution of~\eqref{eq3} near $\tilde u$. Then
$v$ satisfies a linear equation $L_{\hat C}v=0$, where $\hat C$ is an
admissible operator on $N$ in the sense of
Section~\ref{ss:lin-Fredholm}. Moreover, the asymptotic operators of
$\hat C$ at the punctures agree with the asymptotic operators at the
periodic orbits $x_j$ as in (\ref{asy-operator}).
\end{lem}

\begin{proof}
Let us sketch the proof. Denote coordinates on $\dot S$ by $z$ and on
$\R^2$ by $x$. Pick a trivialization $N\cong\dot S\times\R^2$ as
provided by Theorem 4.7 in~\cite{HWZ-Fredholm}. Write the almost
complex structure in this trivialization as
$$
   \tilde J = \left(\begin{matrix} j & \tilde\Delta \\ \Delta & i
   \end{matrix}\right) : T\dot S\times\R^2\to T\dot S\times\R^2,
$$
where the components of $\tilde J$ depend smoothly on $(z,x)$. Since
the zero section $\dot S\times\{0\}$ is $\tilde J$-holomorphic, we
have
$$
   \tilde J(z,0) = \left(\begin{matrix} j(z) & \tilde\Delta(z) \\ 0 &
   i(z) \end{matrix}\right)
$$
for complex structures $j$ on $\dot S$ and $i$ on $\R^2$ and
homomorphisms $\tilde\Delta(z):\R^2\to T_z\dot S$.
According to Section 5 in~\cite{HWZ-Fredholm}, the section $C$
in the bundle $\Hom_\R(N,\A_0)\to\dot S$ defined by
$$
   C(z)h := \Bigl(\frac{\p\Delta}{\p x}(z,0)h\Bigr)\circ j(z,0),
   \qquad h\in\R^2,
$$
is admissible and the corresponding operator
$$
   w\mapsto L_Cw = Tw + i\circ Tw\circ j + Cw
$$
agrees with the normal linearized Cauchy-Riemann operator $D^N$ at
$\tilde u$. In particular, the asymptotic operators of
$C$ at the punctures agree with the asymptotic operators at the
periodic orbits $x_j$.

By hypothesis, the graph $\gr(v)$ of $v$ satisfies the equation
$$
   T\gr(v)+\tilde J\bigl(\gr(v)\bigr)\circ T\gr(v)\circ\hat\jmath = 0
$$
for some complex structure $\hat\jmath$ on $\dot S$. The $\dot
S$-component of this equation yields
$\hat\jmath(z)=j(z)+\Delta(z,v)\circ Tv$; its $\R^2$-component is
$$
    Tv + i(z,v)\circ Tv\circ\hat\jmath + \Delta(z,v)\circ\hat\jmath =
    0.
$$
Define the complex structures $\hat\imath(z) := i\bigl(z,v(z)\bigr)$ on
$\R^2$ and the section $\hat C$ in the bundle $\Hom_\R(N,\A_0)\to\dot S$ by
$$
   \hat C(z)h := \int_0^1\Bigl(\frac{\p\Delta}{\p x}\bigl(z,\tau
   v(z)\bigr)h\Bigr)\circ\hat\jmath\,d\tau,\qquad h\in\R^2.
$$
Then the equation for $v$ can be viewed as a linear
equation as in Section~\ref{ss:lin-Fredholm},
$$
   L_{\hat C}v = Tv + \hat\imath\circ Tv\circ\hat\jmath + \hat Cv = 0.
$$
Note that near a puncture $v(s,t)\to 0$ in $C^\infty$ as
$s\to\infty$,
so $\hat C$ approaches $C$ at the punctures. This implies that $\hat
C$ is admissible with the same asymptotic operators as $C$ and the
lemma follows.
\end{proof}

Finally, let us discuss the effect of exponential weights. By
Theorem~\ref{thm:asymptotics}, $u(s,t)$ approaches $x_j(t)$ normally
with an exponential rate given by an eigenvalue $\lambda_j<0$ of the
asymptotic operator at $x_j$. For weights $\lambda_j<\delta_j<0$ that
are not eigenvalues,
denote by $\MM_\w$ the space of solutions in $\MM$ which
normally approach the $x_j$ at the (positive) punctures with an
exponential rate $\delta_j$ or faster. By construction, the solution
$\tilde u$ belongs to $\MM_\w$. Define the {\em weighted Fredholm index}
$\ind_\w(\tilde u)$ as the index of the linear Fredholm operator
$D:E_\w\to F_\w$ between suitable weighted Sobolov spaces describing
nearby solutions in $\MM_\w$.
An argument as in Section~\ref{ss:lin-Fredholm} shows
(see~\cite{HWZ-Fredholm}, Section 6)
\begin{equation}\label{eq:weighted-full-ind}
   \ind_\w(\tilde u) = \ind(D_\w) = \mu_\w(\tilde u) - 2 + \#\Gamma.
\end{equation}
Similarly, the relations~\eqref{eq:normal-ind}
and~\eqref{eq:ind-normal-ind} carry over to the weighted case.

\subsection{Algebraic invariants}\label{ss:algebraic}
In this section we use the algebraic invariants
from~\cite{HWZ-embedding} to single out a 2-parameter family of
solutions by putting suitable exponential weights. Let $\tilde
u=(a,u):\dot S\to R\times M$ be a (special or generalized) finite
energy sphere with
asymptotic orbits $x_j$. As in the previous section, suppose that
all the punctures are positive and all the asymptotic orbits
are nondegenerate and elliptic with Conley-Zehnder indices $\mu_j =
2\alpha_j + 1$. By Theorem~\ref{thm:asymptotics}, the
solution $u$ approaches $x_j$ from the direction of some
eigenvector $e_j$ of the asymptotic operator. Denote the
winding number of $e_j$ with respect to a trivialization of $u^*\xi$
by $w_j$ and define the {\em asymptotic winding number} of $\tilde u$
by
$$
    \wind_\infty(\tilde u) := \sum w_j.
$$
Since $u$ approaches $x_j$ from the direction of an eigenvector to
a negative eigenvalue at a positive puncture, we have
$$
    w_j\leq\alpha_j.
$$
Let $\tau := \sum(\alpha_j-w_j)$ be the difference between the
actual winding numbers at the punctures and the maximal possible
ones.

Now assume that $\tilde u$ is a special finite energy sphere. Then,
according to~\cite{HWZ-embedding}, the section $\pi\circ Tu$
of the bundle ${\rm Hom}_\C(T\dot S,u^*\xi)$ satisfies a linear
Cauchy-Riemann type equation as in Section~\ref{ss:lin-Fredholm} (this
is not true for generalized finite energy spheres). The $z$-part in
Theorem~\ref{thm:asymptotics} cannot vanish identically because if it
did then $\tilde u$ would be a covering of
the cylinder over $x_j$ and thus have negative punctures (Theorem
6.11. in~\cite{HWZ-embedding}). By the similarity principle, $\pi\circ
Tu$ can only vanish in finitely many points. The 
{\em winding number} $\wind_\pi(\tilde u)$ is then defined as the sum of
the indices of the zeroes of this section. It is a nonnegative integer which
measures how often $u$ is tangent to the Reeb vector field and is
related to the asymptotic winding number by the formula
(see~\cite{HWZ-embedding})
\begin{equation}\label{eq:wind}
    \wind_\pi(\tilde u) = \wind_\infty(\tilde u) - 2 +
    \#\Gamma.
\end{equation}
Combining formulae~\eqref{eq:wind} and~\eqref{eq:full-ind}, we find
\begin{align*}
    2\tau
    &= 2\sum(\alpha_j-w_j)\cr
    &= \sum(2\alpha_j +1)  -
    \#\Gamma - 2\wind_\infty(\tilde u) \cr
    &= \mu(\tilde u) -\#\Gamma - 2\wind_\infty(\tilde u) \cr
    &= \mu(\tilde u) - 2\wind_\pi(\tilde u) - 4 +\#\Gamma \cr
    &= \ind(\tilde u) - 2 - 2\wind_\pi(\tilde u).
\end{align*}
Now pick weights $\delta_j<0$ just above the larger eigenvalue
corresponding to $w_j$. Denote by $\alpha_j^w$ the winding number
corresponding to the largest eigenvalue smaller than
$\delta_j<0$. Then $\alpha_j^w=w_j$ and the sum $\mu_\w(\tilde u)$ of
the corresponding weighted Conley-Zehnder indices satisfies
\begin{eqnarray*}
\mu_w(\tilde{u}) & = & 2\sum\alpha_j^w\,+\,\#\Gamma\\
 & = & 2\sum w_j\,+\,\#\Gamma\\
 & = & -2\tau+2\sum\alpha_j+\#\Gamma\\
 & = & \mu(\tilde{u})-2\tau.
\end{eqnarray*}
The weighted Fredholm
index~\eqref{eq:weighted-ind} becomes
\begin{align*}
   \ind_\w(\tilde u)
   &= \mu_\w(\tilde u) - 2 + \#\Gamma \cr
   &= \ind(\tilde u) - 2\tau \cr
   &= 2 + 2\wind_\pi(\tilde u).
\end{align*}
In particular, if $u$ is embedded and transverse to the Reeb vector
field the winding number $\wind_\pi(\tilde u)$ vanishes and thus
$\ind_\w(\tilde u)=2$.
So we have shown

\begin{lem}\label{lem:weights}
Let $\tilde u=(a,u):\dot S\to\R\times M$ be a special finite energy
sphere having only positive punctures. Suppose that all the asymptotic
orbits are nondegenerate elliptic and that $u$ is embedded and
transverse to the Reeb flow. Then we can introduce exponential
weights at the punctures such that $\tilde u$ belongs to the space
$\MM_\w$ of solutions with these weights, and the Fredholm index of
$\tilde u$ with these weights satisfies $\ind_\w(\tilde
u)=2$.
\end{lem}

\subsection{An implicit function theorem}\label{ss:implicit-function}
Now fix a collection of distinct nondegenerate simple elliptic orbits
$x_j$ and weights $\delta_j<0$. Consider the space $\MM_\w$ of
generalized finite energy spheres with positive punctures asymptotic
to the $x_j$ with exponential decay $\delta_j$ or faster. The
following result describes the local structure of the set $\MM_\w$.

\begin{thm}\label{thm:implicit-function}
Assume that $C=\tilde{u}(\dot S)\in\MM_\w$ is embedded and has
weighted Fredholm index $\ind_\w(\tilde u)=2$. Then neighbouring
solutions in $\MM_\w$ form a smooth 2-dimensional family of mutually
disjoint embedded curves.
\end{thm}

\begin{proof}
The argument is similar to that given in the proof of Theorem 2.7
in~\cite{HWZ2003}. As in Section~\ref{ss:nonlin-Fredholm},
we write neighbouring solutions in $\MM_\w$ as graphs of sections of
the complex normal bundle to $C$. They satisfy a Monge-Ampere type
equation whose linearization $D^N:E_\w\to F_\w$ at the zero section is
a Cauchy-Riemann type operator as in Section~\ref{ss:lin-Fredholm}
between suitable Sobolev spaces with weights $\delta_j$.
Consider a nontrivial element $h$ in the kernel of $D^N$. By
Corollary~\ref{cor:asymptotics}, $h$ approaches zero at the $j$-th
puncture exponentially from a direction $e_j(t)$, where $e_j$ is an
eigenfield of the asymptotic operator at $x_j$. According to
Lemma~\ref{lem:lin-winding}, $h$ has a winding number
$$
   \wind(h) = \sum w_j,
$$
where $w_j$ is the winding number of
$e_j:S^1\rightarrow {\mathbb R}^2\setminus\{0\}$ with respect to a
trivialization of the normal bundle $N$. Let
$$
   \mu_\w^N(x_j)=2\cdot\alpha_\w^N(x_j)+1
$$
be the weighted normal Conley-Zehnder index at $x_j$,
where $\alpha_\w^N(x_j)$ is the maximal winding number of an
eigenfield of the asymptotic operator $A_j$ associated to an
eigenvalue $<\delta_j$. Since $h$ belongs to the Sobolev space with
weights $\delta_j$, we have
$$
   w_j \leq \alpha_\w^N(x_j).
$$
With $\mu_\w^N(\tilde u):=\sum\mu_\w^N(x_j)$ and the
weighted version of formula~\eqref{eq:normal-ind}, this implies
\begin{align*}
   2\wind(h) &= \sum_j 2w_j \cr
   &\leq \sum_j[2\alpha_\w^N(x_j)+1] - \sharp\Gamma \cr
   &= \mu_\w^N(\tilde u) - \sharp\Gamma \cr
   &= \ind_\w(\tilde u)-2 \cr
   &= 0.
\end{align*}
This shows that nontrivial elements in the kernel of $D^N$ are nowhere
vanishing. It follows that the kernel can be at most two-dimensional,
since otherwise we could construct a nontrivial element in the kernel
with a zero. Since $\ind_\w(\tilde u)=2$, we conclude that the
operator $D^N$ is surjective. Thus $\MM_\w$ is a smooth 2-dimensional
manifold near $C$.

It remains to prove that neighbouring elements $C'\neq C$ in $\MM_\w$
do not intersect $C$. As in Section~\ref{ss:nonlin-Fredholm}, describe
$C'$ as the graph of a nonvanishing section $v$ of the normal bundle
to $C$. By Lemma~\ref{lem:nearby-sol}, $v$ satisfies a linear
Cauchy-Riemann type equation $L_{\hat C}v=0$, where the admissible
section $\hat C$ has the same asymptotics as $D^N$. Hence the winding
number of $v$ satisfies
$$
   \wind(v) = \sum w_j.
$$
Now the computation above shows $\wind(v)=0$. Hence $v$ has no zeroes,
which precisely means that its graph does not intersect $C$.
\end{proof}

\begin{defn} For $C\in\MM_\w$ we
denote by $a(C)\in {\mathbb R}$ the minimum of the ${\mathbb
R}$-value of the projection $C\rightarrow {\mathbb R}$.
\end{defn}

Theorem~\ref{thm:implicit-function} has the following immediate
corollary.

\begin{cor}\label{inf}
For $C\in\MM_\w$ as in Theorem~\ref{thm:implicit-function} there
exists a $C'\in\MM_\w$ with
$$
a(C') < a(C).
$$
\end{cor}

\subsection{Intersections}\label{ss:int}
Consider a connected component $\MM_\w^0$ of $\MM_\w$ containing an
embedded solution $C_0$ of index $2$. By
Theorem~\ref{thm:implicit-function} and positivity of intersections
(see~\cite{Floer-Hofer-Salamon}), this implies that $\MM_\w^0$ is a
smooth 2-dimensional manifold and all elements in $\MM_\w^0$ are
embedded. Moreover, nearby distinct elements in $\MM_\w^0$ are
disjoint. The following result shows that any two (not necessarily
nearby) elements are either identical or disjoint.

\begin{prop}\label{prop:int}
Two elements $C,C'$ in $\MM_\w^0$ are either identical or disjoint.
\end{prop}

\begin{proof}
First note that two distinct $C,C'\in\MM^0_\w$ intersect only in
finitely many points. To see this, write $C'$ near the $j$-th puncture
as the graph of a nontrivial section $v_j$ in the normal bundle to
$C$. By Lemma~\ref{lem:nearby-sol} and Lemma~\ref{lem:lin-winding},
$v_j$ has only finitely many zeroes. Thus we have a well-defined
algebraic intersection number $\hbox{int}(C,C')$. Recall
from~\cite{Floer-Hofer-Salamon} that each intersection point
contributes positively to $\hbox{int}(C,C')$ and intersection points
persist under small perturbations. 
 
Now suppose that $C,C'$ in $\MM_\w^0$ are neither identical nor
disjoint, hence $\hbox{int}(C,C')>0$. 
Pick a continuous path 
$(C_{\tau})_{0\le\tau\le 1}$ in $\MM_\w^0$ with $C_0=C$ and
$C_1=C'$. For small $\tau$, we have $\hbox{int}(C,C_{\tau})=0$ by
Theorem~\ref{thm:implicit-function}. 
We define 
\[
\tau_0:=\inf\{0\le\tau\le 1\,|\,\hbox{int}(C,C_\tau)>0\}\in(0,1). 
\]
Since intersections persist under perturbations, we must have
$\hbox{int}(C,C_{\tau_0})=0$. Let $I_{\tau}\subset C$ be
the (finite) set of intersections between $C$ and $C_{\tau}$ on the
surface $C$. Now observe that for every neighborhood $U$ of the set of
punctures on $C$ there is an $\varepsilon>0$ such that   
\[
I_{\tau}\subset U\qquad \text{ for all
}0<\tau <\tau_0+\varepsilon. 
\]
For otherwise we would find a sequence $\tau_n>\tau_0$ with
$\tau_n\to\tau_0$ and intersection points $z_n\in C\cap C_{\tau_n}$
with $z_n\notin U$. But then the $z_n$ would converge to an
intersection point $z\in C\cap C_{\tau_0}$, contradicting
$\hbox{int}(C,C_{\tau_0})=0$.   

As above, write $C_\tau$ for $\tau>0$ near the $j$-th puncture as
the graph of a section $v_j^\tau$ in the normal bundle to $C$ which
approaches zero exponentially from the direction of an eigenfield
$e_j^\tau$ of the asymptotic operator at $x_j$. Denote by
$\wind(e_j^\tau)$ the winding number of $e_j^\tau$ in a trivialization
of the normal bundle, and by $\alpha_j$ the maximal winding
number of an eigenvalue below the weight $\delta_j$. Define an integer
valued function $i(\tau)$ by  
$$
i(\tau) = \hbox{int}(C,C_{\tau})-\sum_{j=1}^{N}
[\hbox{wind}(v^{\tau}_j)-\alpha_j].
$$
If $\tau$ is small the Implicit Function
Theorem~\ref{thm:implicit-function} yields
$\wind(v^{\tau}_j)=\alpha_j$ and 
$\mbox{int}(C,C_{\tau})=0$, hence $i(\tau)=0$. We will show that
$i(\tau)=0$ for all $\tau<\tau_0+\varepsilon$ for some $\eps>0$. Since 
$\sum_{j=1}[\hbox{wind}(v^{\tau}_j)-\alpha_j]\leq 0$, this then implies 
$\hbox{int}(C,C_{\tau})\leq 0$ for all $\tau<\tau_0+\eps$,
contradicting the choice of $\tau_0$. 
 
Fix a sufficiently small neighbourhood $U=\cup_j U_j$ of the set of
punctures on $C$. After trivializing the normal bundle of $C$, the
restriction of each $v^{\tau}_j$ to $\bar U_j$ can be viewed as a complex
valued function on $[0,\infty)\times S^1$ which satisfies a linear
Cauchy Riemann type equation as in
Lemma~\ref{lem:nearby-sol}. Dividing $v^{\tau}_j$ by a suitable smooth
positive function $\beta^{\tau}_j$ and compactifying the infinite half
cylinder to $[0,1]\times S^1$, we obtain functions 
\[
w^{\tau}_j:[0,1]\times S^1\longrightarrow \C
\]
with
\[
\mbox{wind}\bigl(w^{\tau}_j(1,\cdot)\bigr)\,=\,\mbox{wind}(e^{\tau}_j)\,.
\]
These functions are continuous in $\tau$ on the open half-cylinder
$[0,1)\times S^1$. For $\varepsilon$ sufficiently small we have
$I_{\tau}\subset U$ for all $0<\tau<\tau_0+\varepsilon$, hence
$\mbox{wind}\bigl(w^{\tau}_j(0,\cdot)\bigr)$ is well-defined and
independent of $\tau\in(0,\tau_0+\varepsilon)$. As the zeros of
$w^{\tau}_j$ correspond to the intersection points of $C$ and
$C_{\tau}$ and their degrees equal the local intersection number, we 
have
\[
\mbox{wind}\bigl(w^{\tau}_j(1,\cdot)\bigr) - \sum_{\{z:
  w^{\tau}_j(z)=0\}}\mbox{deg}(z)\, =
  \,\mbox{wind}\bigl(w^{\tau}_j(0,\cdot)\bigr). 
\]
We conclude that the left hand side is independent of
$\tau\in(0,\tau_0+\varepsilon)$, and summing over $j$ yields
$i(\tau)=0$ for all $\tau\in(0,\tau_0+\varepsilon)$.
\end{proof}

\section{Proof of the main result}\label{sec:proof}

In this section we prove the main result, Theorem~\ref{thm:main}.

\subsection{Finite energy foliations for Giroux forms}
We start by explaining a result by Abbas~\cite{Abbas} about turning
leaves of an open book decomposition into solutions of
equation~(\ref{eq2}). 

Let us call two solutions $(S,j,\Gamma,\tilde{u})$ and $(S',j',$
$\Gamma',\tilde{u}')$ of (\ref{eq2}) {\em equivalent} if there exists
a biholomorphic map $\phi:(S,j)\rightarrow (S',j')$ mapping $\Gamma$
to $\Gamma'$ (preserving the enumeration) so that
$\tilde{u}'\circ\phi=\tilde{u}$. From now on a solution of our
differential equation is an equivalence class
$[S,j,\Gamma,\tilde{u}]$. Note that we have a natural ${\mathbb
R}$-action on the solution set by associating to $c\in{\mathbb R}$
and $[S,j,\Gamma,\tilde{u}]$ the new solution
$$
c+[S,j,\Gamma,\tilde{u}] := [S,j,\Gamma,(a+c,u)],\qquad \tu=(a,u).
$$
Given $[S,j,\Gamma,\tilde{u}]$, we denote by
$C$ the image of $\tilde{u}$. Since all the maps $\tilde{u}$ of
interest to us will be somewhere injective one can show that
knowing $C$ we can reconstruct the underlying equivalence class
$[S,j,\Gamma,\tilde{u}]$. A crucial concept for our discussion is
the notion of a finite energy foliation ${\mathcal F}$.

\begin{defn}
A smooth foliation ${\mathcal F}$ of ${\mathbb R}\times M$ is
called a {\em finite energy foliation} if every leaf $F$ is the image
of an embedded solution $[S,j,\Gamma,\tilde{u}]$ of~(\ref{eq2}),
$$
F=\tilde{u}(\dot{S}),
$$
and with every leaf $F\in{\mathcal F}$ also $c+F\in{\mathcal F}$
for every $c\in{\mathbb R}$, i.e., the foliation is ${\mathbb
R}$-invariant.
\end{defn}

Finite energy foliations are known to be a useful tool in studying
the dynamics of Reeb vector fields as well as topological
applications, see
\cite{HWZ2003,HWZ-strictly-convex,HWZ-characterization}. The following
theorem is proved by Abbas in~\cite{Abbas}.

\begin{thm}\label{thm:Abbas}
Assume that $M$ is a closed three-manifold equipped with a planar
contact structure $\xi$ and a Giroux form $\lambda_{Giroux}$ with
nondegenerate closed Reeb orbits. Let $J$ be any compatible
complex multiplication on $\xi$. Then there exist an open book
decomposition $(L,\pr)$ for $\xi$ with Giroux form $\lambda_{Giroux}$
and a finite energy foliation ${\mathcal F}$ of ${\mathbb R}\times M$
with the following properties.
\begin{itemize}
\item The cylinders over the binding orbits in $L$ are leaves of $\FF$,
  called the {\em trivial leaves}.
\item Every nontrivial leaf is the image of a finite energy sphere
  with only positive punctures. These punctures are in 1-1
correspondence with, and asymptotic to, the binding orbits.
\item The projection to $M$ of any nontrivial leaf can be compactified
  (by the binding orbits) 
to a page of the open book decomposition $(L,\pr)$.
\end{itemize}
\end{thm}

\begin{rem} 
(1) The open book decomposition $(L,\pr)$ in the theorem may differ
from the planar open book decomposition we started with.
 
(2) Theorem~\ref{thm:Abbas} can be proved along the following lines,
using the compactness theorem for symplectic field
theory~\cite{BEHWZ}.   
The first step consists of modifying a leaf $u_0$ of the given planar
open book decomposition $(L_0,\pr_0)$ near its punctures so that there
is a suitable function $a_0$ such that $\tilde{u}_0=(a_0,u_0)$ solves
the differential equation near the punctures. This is achieved by
choosing $J:\xi\rightarrow\xi$  
near the binding $L_0$ and the complex structure $j$ near the punctures 
in a very special way so that solutions can be written down explicitly. 
We then look for a global solution $\tilde{u}=(a,u)$ to the equation 
$u^{\ast}\lambda\circ j=da$ of the form $u=\phi_f(u_0)$, where $f$ is a 
suitable real valued function on the closed surface $S$ (a sphere in our 
case) and where $\phi_t$ denotes the flow of the Reeb vector field. This 
amounts to solving an inhomogeneous Cauchy Riemann equation for the 
function $a-a_0+if$ on the sphere, which is possible because on the
sphere $\overline{\partial}$ is surjective. The first part
of~(\ref{eq2}) involving $\pi\circ Tu$ can then be used to define a
$z$--dependent complex structure $J^+$ on $\xi$ so that
\[
\pi Tu(z)\circ j\,=\,J^+(z,u(z))\circ \pi Tu(z).
\]
A cobordism argument similar to the one in Section 3 of this paper 
can then be used to deform the parameter dependent complex structure 
$J^+$ into one which does not depend explicitly on $z$, say $J^-$. We 
pick a complex structure $\tilde{J}=\tilde{J}(z,a,u)$, $(z,a,u)\in 
S\times{\mathbb R}\times M$, on ${\mathbb R}\times M$ such that 
$\tilde{J}\equiv \tilde{J}^+$ for $a\ge 1$ and $\tilde{J}\equiv 
\tilde{J}^-$ for $a\le 0$, and we study the corresponding
PDE~(\ref{eq3}). There is an implicit function theorem and the
compactness result~\cite{BEHWZ} can be applied. Assume that
$\tilde{u}_k=(a_k,u_k)$ is a sequence of solutions such that
\[
\inf a_k\rightarrow r\in{\mathbb R}.
\]
Although there is no statement corresponding to Theorem 3.4 in this 
paper, a solution where the infimum equals $r$ can still be found (the 
part in the cobordism of the broken punctured holomorphic curve in the
limit). An argument similar to the one in section 3.4 of this paper
produces a finite energy solution $\tilde{u}=(a,u)$ to the PDE in the
negative part $({\mathbb R}\times M,\tilde{J}^-)$ with only positive
punctures 
such that $u$ is an embedding transverse to the Reeb vector field. The
collection $L$ of positive punctures of $u$ may differ from the
binding $L_0$ of the original open book decomposition. It is
then shown that there is a compact 1--dimensional family of such
solutions which form an open book decomposition with binding $L$.
\end{rem}

We will refer to the nontrivial leaves in Theorem~\ref{thm:Abbas} as
{\em Abbas solutions}.

\subsection{A cobordism}\label{ss:cobordism}
Suppose now that $\xi$ is supported by a planar open book
decomposition. Let $\lambda_{Giroux}$ be an associated Giroux form
with nondegenerate elliptic binding orbits. We are interested in the
Reeb flow of a different contact form $\lambda$ defining
$\xi$. Multiplying the Giroux form by some positive constant, we may
assume that
$$
   \lambda_\Giroux = f^+\cdot\lambda
$$
for a function $f^+>1$ on $M$. Pick $R>0$ and $\tilde J$ as in the
beginning of Section~\ref{sec:finite-energy-spheres} (with
$\lambda^+=\lambda_\Giroux$ and $\lambda^-=\lambda$) and consider the
PDE~\eqref{eq3} for generalized finite energy spheres. Observe that
any Abbas solution is, after
translating it by a sufficiently large positive constant, a solution
of~\eqref{eq3}. Denote by ${\mathcal A}$ the collection of
all images of Abbas solutions which are contained in
$[R,\infty)\times M$. Of course, any two such solutions are either
disjoint or identical, and the space $\A$ is connected.

Let $A=\tilde u(\dot S)$ be an Abbas solution. Pick weights $\delta_j$
as in Lemma~\ref{lem:weights} so that $A$ has weighted Fredholm index
$\ind_\w(A)=2$. Denote by $\MM_\w$ the space of solutions of
equation~\eqref{eq3} with positive punctures asymptotic to the
binding orbits of $\lambda_\Giroux$ and with weights $\delta_j$.
Note that since all Abbas solutions have the same winding numbers at
the punctures, we have $\A\subset\MM_\w$. Let $\MM_\w^0$ be the
connected component of $\MM_\w$ containing $\A$.

\subsection{A compactness statement}
Assume for the moment that $\lambda$ is nondegenerate. Then we have
the following compactness result for $\MM_\w^0$.

\begin{thm}\label{thm:comp}
Assume that $\lambda$ is nondegenerate.
Let $C_k=\tilde u_k(\dot S)$ be a sequence in $\MM_\w^0$ so that
$a(C_k)\rightarrow r\in {\mathbb R}$. Then, after passing to a
subsequence, there exists an element
$C\in\MM_\w^0$ so that for suitable parametrizations
$\tilde{u}_k\rightarrow \tilde{u}$ in $C^{\infty}_{loc}$. Moreover,
$a(C)=\lim_{k\rightarrow\infty} a(C_k)$.
\end{thm}

\begin{proof}
We apply the compactness theorem for symplectic field
theory~\cite{BEHWZ}. After passing to a subsequence, the $C_k$
converge to a broken punctured holomorphic curve of type
$(k^-|k^0|k^+)$. This means that the limit curve has
$k^+\geq 0$ components in the symplectization $(\R\times M,\tilde
J^+)$ of the positive end, $k^-\geq 0$ components in the
symplectization $(\R\times M,\tilde J^-)$ of the negative end, and
$k^0\in\{0,1\}$ components in the cobordism $(\R\times M,\tilde J)$.
From $a(C_k)\rightarrow r\in {\mathbb R}$ we conclude that $k^-=0$ and
$k^0=1$, so the limit curve is of type $(0|1|k^+)$. If $k^+=0$ the
assertion of the proposition follows from the definition of
convergence in~\cite{BEHWZ}. Thus suppose that $k^+\geq 1$. Then the
top layer is a (not necessarily connected) curve solving the
homogeneous $\tilde{J}^+$-problem. By the stability requirement
in~\cite{BEHWZ}, this layer contains at least one component $\hat C$
that is not a cylinder over a closed Reeb orbit. Let $\hat u$ be a
parametrization of $\hat{C}$. Note that $\hat u$ must be somewhere
injective since the positive asymptotic limits are simply-covered
(they are binding orbits of the open book decomposition). Therefore,
$\hat u$ is not a branched covering of a cylinder over a closed Reeb
orbit, and hence $\hat C$ has to intersect a nontrivial leaf of the
Abbas foliation. On 
the other hand, $\hat C$ cannot be identical to such a leaf since it
has at least one negative puncture. From the definition of convergence
and positivity of intersections it follows that there is a sequence of
Abbas solutions $A_k\in\A$ such that $C_k\cap A_k\neq \emptyset$ and
$C_k\neq A_k$ for large $k$. Since $\A\subset\MM_\w^0$, this
contradicts Proposition~\ref{prop:int}.
\end{proof}

\subsection{Conclusion}
The planar Weinstein conjecture is now proved as follows. We keep
assuming that $\lambda$ is nondegenerate. First note that
$$
   \inf_{C\in\MM_\w^0} a(C) =-\infty.
$$
Indeed, arguing indirectly, suppose the left-hand side defines
a real number $r$. Take a sequence $C_k\in\MM_\w^0$ with
$a(C_k)\rightarrow r$. By Theorem~\ref{thm:comp}, after taking a
subsequence, we find a $C\in\MM_\w^0$ with
$$
   a(C) = \lim_{k\to\infty}a(C_k) = \inf_{C'\in\MM_\w^0} a(C').
$$
By Corollary~\ref{inf}, there exists a $C'\in\MM_\w^0$ with
$a(C')< a(C)$, giving a contradiction. The proof is now completed by
taking a sequence $C_k\in\MM_\w^0$ with
$$
   a(C_k)\rightarrow-\infty.
$$
We apply again the compactness theorem for symplectic field
theory~\cite{BEHWZ}. After passing to a subsequence, the $C_k$
converge to a broken punctured holomorphic curve of type
$(k_-|k^0|k_+)$. From $a(C_k)\rightarrow -\infty$ we conclude $k_-\geq
1$ (and consequently $k^0=1$). By the definition of convergence
in~\cite{BEHWZ}, the lowest layer of the limit curve contains a
non-constant (special) finite energy sphere $\hat C$ for
$(M,\lambda,J^-)$ having only positive punctures. At the punctures,
$\hat C$ is asymptotic to periodic Reeb orbits $x_j$ of the contact
form $\lambda$. By construction, their homology classes $[x_j]$
satisfy $\sum [x_j]=0$, where we sum over all the punctures. This
proves Theorem~\ref{thm:main} in the case that $\lambda$ is
nondegenerate.

If $\lambda$ is degenerate we can take a sequence
$f^{(k)}:M\rightarrow (0,\infty)$ of smooth functions converging in
$C^{\infty}$ to the constant function $f(x)\equiv 1$ so that
the contact forms $f^{(k)}\lambda$ are nondegenerate. By the result in
the nondegenerate case, we find for every $k$ a finite set of periodic
Reeb orbits $x_j^{(k)}$ for $f^{(k)}\lambda$ whose homology classes
sum up to zero. By the proof of the compactness theorem
in~\cite{BEHWZ}, the number of orbits for each $k$ is bounded by a
constant independent of $k$. So after passing to a subsequence, we may
assume that their number is constant. Then by the Arzela-Ascoli
theorem, after passing to a subsequence, the $x_j^{(k)}$ converge in
$C^\infty$ as $k\to\infty$ to periodic Reeb orbits for the contact
form $\lambda$. Clearly, the homology classes of the $x_j$ sum up to
zero. This concludes the proof of Theorem~\ref{thm:main}.

\end{document}